\documentclass[11pt]{amsart}
\usepackage{amsaddr}
\usepackage{geometry}                
\usepackage{graphicx}
\usepackage{amssymb}
\usepackage{epstopdf}
\usepackage{amssymb}
\usepackage{amsmath}
\usepackage{blkarray} 
\usepackage{comment}
\usepackage{latexsym}
\usepackage{amsthm} 
\usepackage{eucal} 
\usepackage{eufrak}
\usepackage{float}
\usepackage[rightcaption]{sidecap}
\usepackage{graphicx}
\graphicspath{ {images/} }
\usepackage[utf8]{inputenc}
\usepackage[main=english, ukrainian]{babel}
\usepackage{xcolor}
\usepackage{mathrsfs}
\usepackage{array}
\usepackage{framed}
\usepackage[all]{xy}
\usepackage{sidecap}
\usepackage{wrapfig}
\usepackage{multicol}
\usepackage{lscape}
\usepackage{MnSymbol}
\usepackage{mathtools}
\usepackage{cite}
\usepackage{appendix}
\usepackage{chngcntr}
\usepackage{etoolbox}
\usepackage{lipsum}
\usepackage{rotating}
\usepackage[thinlines]{easytable}
\usepackage{hyperref}
\usepackage{empheq}
\usepackage{enumitem}
\usepackage{lineno}

\DeclareMathAlphabet{\mathpzc}{OT1}{pzc}{m}{it}

\theoremstyle{plain} 
\newtheorem{thm}{Theorem}[section] 
\newtheorem{cor}[thm]{Corollary} 
\newtheorem{lem}[thm]{Lemma} 
\newtheorem{prop}[thm]{Proposition} 

\theoremstyle{definition} 
\newtheorem{defn}{Definition}

\theoremstyle{remark} 

\newtheorem{rmk}{Remark}

\title{Good and bad children in metabolic networks}
\author{Nicola Vassena}
\address{Freie Universit{\"a}t Berlin}
\email{nicola.vassena@fu-berlin.de}
\date{\today}

\begin{document}

\providecommand{\keywords}[1]
{
  \small	
  \textbf{\textit{Keywords---}} #1
}

\maketitle

\tableofcontents

\begin{abstract}
Equilibrium bifurcations arise from sign changes of Jacobian determinants, as parameters are varied. Therefore we address the Jacobian determinant for metabolic networks with general reaction kinetics. \\
Our approach is based on the concept of Child Selections: each (mother) metabolite is mapped, injectively, to one of those (child) reactions that it drives as an input.\\
Our analysis distinguishes reaction network Jacobians with constant sign from the bifurcation case, where that sign depends on specific reaction rates.\\
In particular, we distinguish ``good" Child Selections, which do not affect the sign, from more interesting and mischievous ``bad" children, which gang up towards sign changes, instability, and bifurcation.\\

\smallskip

\end{abstract}

\section{Introduction} \label{intro}

One of the major obstacles in the study of dynamical systems arising from chemical reaction networks is the lack of precise quantitative knowledge of the parameters involved. For this reason, a natural approach is to connect dynamical properties with the network structure, only. In fact, this idea has a long tradition, in chemical reactions settings. Pioneering works by Horn\&Jackson \cite{HFJ72} and Feinberg \cite{Fei87, Fei95} have been directed towards sufficient network conditions to assert existence and uniqueness of a stable equilibrium. The concept of \emph{deficiency} has been developed towards this purpose. Further network features such as \emph{injectivity} \cite{GaNi65, CraFei06, Ba-07, BaCra10} and \emph{concordance} \cite{ShiFei12, ShiFei13} emerged in this context to forbid multistationarity. In particular, saddle-node bifurcations are excluded.\\
In another direction, in 1981 Thomas \cite{Thom81} conjectured a necessary network condition for multistationarity. In particular, certain \emph{positive loops} in the network were identified as necessary network motives to sustain multistationarity. The conjecture was then fully proven in mathematical settings by Soul\'{e} \cite{Soule2003}. Refined conditions for restricted classes of networks have been further studied till present days. See \cite{KST07, WF13, Timo20}, among many others.\\ 
An excellent compendium of the various approaches to multistationarity questions can be found in \cite{MR07}, where the authors also presented and clarified the important related work of Ivanova \cite{Iv79,Iv792, VIv87}, independently developed in the `70s and unfortunately not easily accessible to english readers.\\
Focusing on metabolic networks, the novelty of our present contribution is in two regards. Firstly, we fully characterize on a graphical level the sign of the  determinant of the Jacobian matrix. Secondly, we provide a recipe to detect sign changes of the Jacobian determinant in quite general examples. This implies the identification of parameter areas for equilibria bifurcations, leading to stability changes and multistationarity.\\

\begin{figure} 
\begin{center}
\includegraphics[scale=0.7]{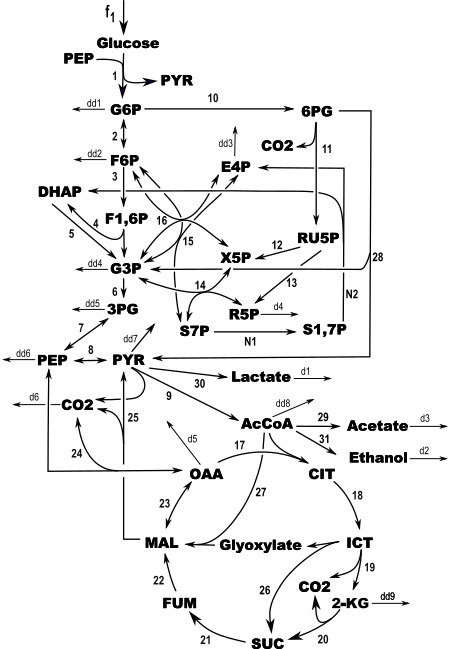}
\caption{The network represents the \emph{central carbon metabolism} of \textit{Escherichia coli}. This figure has been taken from \cite{BF18} and the graphical representation is courtesy of Anna Karnauhova. The inflow feed reaction is named $f_1$. Outflow exit reactions are labeled $d1-d6$ and $dd1-dd9$. Here, for image simplicity, a reversible arrow reaction $m\longleftrightarrow \tilde{m}$ encodes two different opposite reactions. Metabolites $PEP$, $PYR$ and $CO2$ have been graphically repeated, for sake of clarity of the picture.}
\label{anna}
\end{center}
\end{figure}

Metabolism is central to life. A metabolic process is a sequence of chemical reactions designed to transform nutrients into energy. A prominent example of a metabolic network is the \emph{central carbon metabolism} of \textit{Escherichia coli}, figure \ref{anna}. The network depicts the metabolic transformation of glucose, which allows the production of energy. Furthermore, this network is widely used as a general model for cellular respiration, due to the fact that \textit{E. coli} are genetically well-known and easy to treat. We will concentrate on this network in Section \ref{TCAC}. However, this reference example helps us at first to focus on some mathematical features of metabolic networks.\\

In particular, we underline (\hypertarget{star}{$\star$}):

\begin{enumerate} 
\item The stoichiometric coefficients in metabolic networks are mostly \{1,0\}. In the \textit{E. coli} network in figure \ref{anna} they are \emph{only} \{1,0\}.
\item The number of metabolites involved in a reaction is considerably low compared to the size of the network. The \textit{E. coli} network in figure \ref{anna} consists of 30 metabolites interacting through 58 reactions. The maximum number of metabolites involved in a reaction is four (e.g. reaction 14), and many reactions involve only two metabolites (e.g. reaction 10).
\item The reaction rate functions (kinetics) describe the mathematical form of the reactions. Metabolic networks, such as the one in figure \ref{anna}, represent only the metabolic transformations without considering the enzymes network. That is, enzymes and other secondary chemicals do not appear explicitly in the network. The enzymes are usually taken in account by using enzymatic kinetics (e.g. Michaelis-Menten) instead of elementary kinetics (e.g. mass action).
\end{enumerate}
Using the above outline as a guide, we can begin to discuss the details of our approach.\\

We consider general metabolic chemical reaction networks $\mathbf{\Gamma}$ with $M$ metabolites and $N$ reactions. For notation, we use labels $A, B, C, D, ...$ for metabolites and $1, 2, 3, ...$ for reactions. We call $\textbf{M}$ the set of metabolites and $\textbf{E}$ the set of reactions, such that $|\mathbf{M}|=M$ and $|\mathbf{E}|=N$. We use the small letter $m \in \textbf{M}$ for a generic metabolite and the small letter $j\in \mathbf{E}$ for a generic reaction.\\

A chemical reaction $j$ is represented as
\begin{equation}
 j: \quad s^{j}_1m_1+...+s^{j}_Mm_M \underset{j}{\longrightarrow} \bar{s}^{j}_1m_1+...+\bar{s}^{j}_Mm_M,
\end{equation}
with nonnegative stoichiometric coefficients $s^{j},\bar{s}^j \in \mathbb{R}$. In a metabolic context, we repeat, these coefficients are mostly 0 or 1.\\
A metabolite $m$ is called an \emph{input} or a \emph{reactant} of the reaction $j$, if $s^{j}_m \neq 0$. Respectively, $m$ is called  an \emph{output} or a \emph{product} of the reaction $j$, if $\bar{s}^{j}_m\neq 0$. We say, conversely, that a reaction $j$ is \emph{outgoing} from the metabolite $m$ if $m$ is an input of reaction $j$. We say that a reaction $j$ is an \emph{ingoing} reaction of the metabolite $m$ if $m$ is an output of $j$.\\
Metabolic systems are intrinsically open systems, that is, they exchange chemicals with the outside environment by feed and exit reactions. Within our settings, the constant \emph{feed reactions}, or \emph{inflows}, are reactions with no inputs ($s^j=0$) and the \emph{exit reactions}, or \emph{outflows}, are reactions with no outputs ($\bar{s}^j=0$).\\

Graphically, we represent a reaction
\begin{equation}\label{firstpiceq}
j: \quad  A+2B \underset{j}{\longrightarrow} C,
\end{equation}
where we have omitted stoichiometrically zero terms, as follows.
\begin{equation}\label{firstpic}
\end{equation}
\begin{center} 
\vspace{-0.7cm}
\includegraphics[scale=0.3]{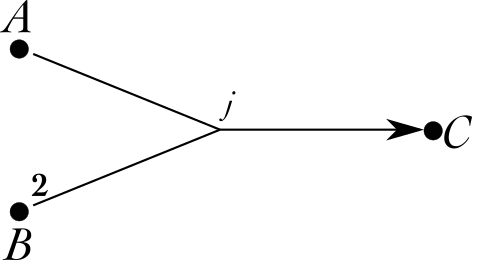}
\end{center}

In \eqref{firstpic}, the arrow orientation is inherited from \eqref{firstpiceq}. The stoichiometric coefficient 2 of metabolite $B$ is indicated as a weight in the lower tail of the directed arrow $j$, and stoichiometric coefficients 1 are omitted, as well as non-participating other reactants. In particular, this graphical representation considers the metabolites as vertices and the reactions as arrows of the network, and it is one natural representation widely used in chemistry, biology, and mathematics.\\

Explicit autocatalytic reactions $j$ are defined as reactions for which a metabolite $m$ is both an input and an output of the reaction. In symbols,  $s^j_m,\bar{s}^j_m \neq 0$, for at least one metabolite $m$. \emph{Throughout this paper, we exclude explicit autocatalytic reactions}. In particular, self-loops are not allowed in the graphical representation of the network.\\
This assumption is only made for mathematical simplicity. Explicit autocatalysis can be treated via intermediaries. For example, the reaction 
\begin{equation}  \label{autocat}
j: \quad m \underset{j}{\longrightarrow} 2m
\end{equation}
is explicit autocatalytic. However, the same chemical process can be modeled without explicit autocatalytic reactions as
\begin{equation}  
m \underset{j_1}{\longrightarrow} m_0 \underset{j_2}{\longrightarrow} 2m,
\end{equation}
with $m_0$ as an intermediary.\\

To construct the $M\times N$ stoichiometric matrix $S$, let us consider any reaction $j$. We associate to any stoichiometric coefficient $s^j_m$ of an \emph{input} metabolite $m$ of the reaction $j$ a \emph{negative} stoichiometric entry of the stoichiometric matrix S, that is: 
\begin{equation}
S_{mj}:=-s^j_m\,,\quad \text{ for $m$ input of $j$}.
\end{equation}
Conversely, we associate to any stoichiometric coefficient $\bar{s}^j_{m}$ of an \emph{output} metabolite $m$ of the reaction $j$ a \emph{positive} stoichiometric entry of S, that is:
\begin{equation}
S_{mj}:=\bar{s}^j_{m}\,,\;\quad \text{ for $m$ output of $j$}.
\end{equation}
For example, a \emph{monomolecular reaction} $j$ is a reaction which possesses as input one single metabolite $m_1$ and as output one single metabolite $m_2$,
\begin{equation}
 m_1 \underset{j}{\longrightarrow} m_2.
\end{equation}
Such a reaction translates into the $j^{th}$ column of the stoichiometric matrix $S$ as
\begin{equation}
S^j=
\begin{blockarray}{cc}
 & j \\
\begin{block}{c(c)}
  m_1 & -1\\
  m_2 & 1\\
  m_3 & 0\\
  ... & ...\\
  m_M & 0\\
\end{block}
\end{blockarray}\;.
\end{equation}
Here we have indicated the rows by $m_1, ... , m_M$, explicitly. In particular, with this construction we model a reversible reaction
\begin{equation}
j: \quad A+B \underset{j}{\longleftrightarrow} C
\end{equation}
simply as two distinct irreversible reactions
\begin{equation}
j_1: \quad A+B \underset{j_1}{\longrightarrow} C\quad \text{ and }\quad j_2: \quad C \underset{j_2}{\longrightarrow} A+B.
\end{equation}
Columns associated with feed reactions contain only positive entries and the ones associated with exit reactions contain only negative entries. All other columns contain both positive and negative entries. In the \textit{E. coli} picture of figure \ref{anna}, only for image simplicity, a reversible arrow reaction $m\longleftrightarrow \tilde{m}$ encodes two different opposite reactions.\\
Stoichiometric matrices of metabolic networks, due to our observations (\hyperlink{star}{$\star$}), are then \emph{sparse} matrices with most entries being \{-1, 0, 1\}.\\

Let $x_m(t)$ be the time evolution of the concentration of the metabolite $m$. The isothermal dynamics of the vector $x \in \mathbb{R}^M$ of the concentrations is described by the system of differential equations
\begin{equation} \label{eq}
\dot{x}=f(x):=S \mathbf{r}(x).
\end{equation}
The $M \times N$ matrix $S$ is the stoichiometric matrix constructed above. The $N$-dimensional vector $\mathbf{r}(x)$ represents the reaction rates as functions of $x$: the \emph{kinetics} of the system. The feed reactions are represented by constant functions; i.e.:
\begin{equation}
r_{j_f}(x)\equiv K_{j_f},
\end{equation}
for a feed reaction $j_f$. 
\emph{Throughout this paper, we pose the following assumptions on the reaction rates }$\mathbf{r}(x)$:

\begin{enumerate}
\item We assume the reaction rates $r_j(x)$ to depend only on those concentrations $x_m$ such that the metabolite $m$ is an input metabolite of reaction $j$. In particular,
\begin{equation*}
\frac{\partial r_j(x)}{\partial x_m}\equiv 0,\quad \text{ unless $m$ is an input of $j$}.
\end{equation*}
Moreover, we use the notation $r_{jm}$ for the nonzero partial derivatives, i.e.,
\begin{equation*}
r_{jm}:= \frac{\partial r_j(x)}{\partial x_m} \neq 0
\end{equation*}
if $m$ is an input of reaction $j$.
\item We consider strictly positive monotone reaction rate functions $r_j(x) \in C^1$, for every $j=1,...,N$:
\begin{equation*}
r_j(x)>0 \quad\text{ for }\quad x>0,
\end{equation*}
and, for the nonzero partial derivatives $r_{jm}$, strictly positive slopes
\begin{equation*}
 r_{jm} > 0.
\end{equation*} 
This monotonicity restriction is indeed satisfied for most, but not all, chemical reaction schemes. 


\end{enumerate}

With these assumptions, all required information of the network is encoded in the stoichiometric matrix $S$, only. In particular, we do not specify here the mathematical form of the kinetics, remaining in wide generalities.\\


Great effort in the mathematical community has been spent in finding network characterizations of the existence and the uniqueness of equilibrium solutions of \eqref{eq}. For extensive reference, see the comprehensive book by Feinberg \cite{Fei19}.\\
With our approach, we do not address this question at all. In fact, \emph{throughout this paper, we assume the existence of a dynamical equilibrium} $x^*$ that solves 
\begin{equation} \label{eqss}
0=f(x^*):=S \mathbf{r}(x^*).
\end{equation}
The assumption of the existence of a dynamical equilibrium is not smoothly untroubled. In particular, linear constraints have been implicitly imposed on the reaction rates $\mathbf{r}$, because of \eqref{eqss}. Note that these constraints do not necessarily fix the precise value of an equilibrium $x^*$, and can be considered posed a priori, so that the existence of the equilibrium is an assumption on the reaction rates $\mathbf{r}$, only. Here, our analysis is based entirely on the derivatives $r_{jm}$ of the reaction rates and we do not want to be concerned by the {equilibrium} constraints \eqref{eqss}. To avoid this, we must assume a certain independence of the derivatives $r_{jm}$ from the reaction rates themselves, at the equilibrium. In particular, locally at the fixed equilibrium, we require the possibility of freely choosing the value of any $r_{jm}$ independently from each other and from the constraints $S\mathbf{r}=0$. In this sense, the partial derivatives $r_{jm}$ can be considered \emph{positive free parameters}. This requires a certain mathematical complexity of the reaction rates $r_j$. In fact: too mathematically `simple' kinetics fail to satisfy this assumption.\\ 

As an example, for polynomial mass action kinetics, the value of $r_j(x)$ and $r_{jm}(x)$ are related, a priori, at any value $x$, and for any $j$ and $m$. In particular, the theory developed here does not {fully} apply to mass action kinetics. In contrast, Michaelis-Menten kinetics, more suited for metabolic networks, satisfies our independence assumption. {We present here an exemplification of this fact. The complete mathematical argument can be found in \cite{F19}, which is a contribution by Fiedler that shares this premise. Let us consider the following reaction, whose single input is a metabolite $m$:
\begin{equation}
m \quad\mathrel{\mathop{\longrightarrow}^{\mathrm{1}}} \quad ...
\end{equation}
The product of the reaction is irrelevant for this discussion, and it is omitted. The rate of the reaction 1, according to the law of \emph{mass action}, reads
\begin{equation}
r_1(x_m)=k_1x_m.
\end{equation}
Here, $k_1$ is a positive coefficient. The derivative of $r_1(x_m)$ with respect to $x_m$ is then
\begin{equation}
r_{1m}(x_m)=\frac{\partial r_{1}}{\partial x_m}(x_m)=k_1.
\end{equation}
In particular, at \emph{any} fixed $\bar{x}_m$, the value $r_1(\bar{x}_m)$ uniquely determines the derivative $r_{1m}(\bar{x}_m)$. Indeed,
\begin{equation}
r_{1m}(\bar{x}_m)=\frac{r_1(\bar{x}_m)}{\bar{x}_m}.
\end{equation}
For general nonlinearities, of course, we may choose $r_1(\bar{x}_m)$ and $r_{1m}(\bar{x}_m$) independently of each other. In the special case of \emph{Michaelis-Menten} kinetics, for example, the rate of reaction 1 is
\begin{equation}
r_1(x_m)=\frac{k_1x_m}{1+a_1x_m},
\end{equation}
where both $k_1$ and $a_1$ are positive parameters, and its derivative reads
\begin{equation}
r_{1m}(x_m)=\frac{k_1}{(1+a_1x_m)^2}.
\end{equation}
This implies that at \emph{any} fixed $\bar{x}_m$, the value $r_{1m}(\bar{x}_m)$ can be chosen as small as needed, independently from $r_1(\bar{x}_m)$:
\begin{equation}
\frac{r_{1m}(\bar{x}_m)}{r_1(\bar{x}_m)}=\frac{1}{(1+a_1\bar{x}_m)\bar{x}_m} \in (0,1)\cdot \frac{1}{\bar{x}_m}.
\end{equation}
In particular, for large $a_1 \rightarrow \infty$ we get small $r_{1m}\rightarrow 0$, and for small $a_1 \rightarrow 0$ we get $r_{1m} \rightarrow \frac{r_1(\bar{x}_m)}{\bar{x}_m}$. Having fixed the value $\bar{x}_m$ requires then to pick $k_1:=\frac{(1+a_1\bar{x}_m)r_1(\bar{x}_m)}{\bar{x}_m}$, for any chosen $a_1$. Note that for the limit value $a_1=0$ we recover the mass action case.\\
Let us concentrate now on the relevant case in which the fixed $\bar{x}$ is an equilibrium. In the mass action case, it is not possible to separate the question of the existence of the equilibrium (related to the reaction rates $\mathbf{r}(\bar{x}$)) from the question of its stability (related to the derivatives $r_{jm}(\bar{x})$), because the former determines the latter as shown in the analysis above. On the contrary, for more general kinetics as Michaelis-Menten, the two questions can be addressed separately and independently.\\}

{The Jacobian matrix of a dynamical system plays a central role in the stability analysis of equilibria. The sign of the eigenvalues is an indication of stability; therefore, a change in sign of the determinant hints at a change in stability and bifurcation phenomena.} Let $G$ be the Jacobian matrix of the equilibrium system \ref{eqss}, that is:
\begin{equation}
G:= f_x, \quad \quad \quad \text{with entries}  \quad \quad \quad G_{kh}= \frac{\partial f_k}  {\partial x_h}.
\end{equation}

Note that the entries $G_{ij}$ are linear forms in the variables $r_{jm}$.\\ 

The leading question of this paper is the following:
\begin{equation} \label{question}
\textit{When is $\operatorname{det}G$ of fixed sign?}
\end{equation}
That is: 
\begin{center}
\textit{When - for \emph{any} choice of positive parameters $r_{jm}$ - does the determinant carry the same sign?}
\end{center}

In this paper, we fully characterize, on a graphical level, the answer to the question \eqref{question}. In particular, in Section \ref{CBA}, we introduce \emph{Child Selections} and we use the Cauchy-Binet formula to expand the Jacobian determinant in a polynomial, in which each monomial summand is associated to a Child Selection. Depending on the sign of the coefficients of these monomials, each Child Selection is abstractly classified in \emph{good} or \emph{bad}. The coexistence of a good and a bad Child Selection characterizes the condition of \emph{indeterminate sign} Jacobian. In Section \ref{MT1}, the main Theorem \ref{MainThm1} abstractly characterizes whether a given Child Selection is good or bad and Section \ref{IR1} translates this abstract condition into a graphical network condition. Consequently, Section \ref{SNB} uses the developed concepts to find a bifurcation parameter responsible for a change of sign in the Jacobian determinant, with possible consequent bifurcation phenomena. Section \ref{TCAC} contains an example of an application for the central metabolism of {\textit{E. coli}}. Section \ref{discussion} and \ref{proofs} conclude with the discussion and proofs.\\

\textbf{Acknowledgments.} We thank the anonymous referees for their helpful and valuable comments. We are very grateful to Bernold Fiedler for his enduring support. Our thanks go also to the present and past group of Nonlinear Dynamics of Freie Universit{\"a}t Berlin. Anna Karnauhova with everlasting gentleness allowed us to use her masterful picture of the central metabolism of {\textit{E. coli}}. Meghan Kane greatly helped in fixing non-native english typos. This work was funded by the Deutsche Forschungsgemeinschaft (DFG, German Research Foundation)-project number 163436311-SFB 910, and supported by the Berlin Mathematical School.

\section{Cauchy-Binet analysis via Child Selections} \label{CBA}

For the metabolic chemical reaction system \ref{eq}
\begin{equation*} 
\dot{x}=f(x):=S \mathbf{r}(x),
\end{equation*}
the Jacobian matrix reads
\begin{equation} \label{eqdetfirst}
G=f_x=SR.
\end{equation}
The \emph{reactivity matrix} $R$ of the partial derivatives is an $N \times M$ matrix, whose entries $R_{jm}$ are given by:
\begin{equation} \label{Req}
R_{jm}:=\frac{\partial}{\partial x_m} r_j(x) =
\begin{cases}
r_{jm}\quad \text{if}\quad\frac{\partial r_j(x)}{\partial x_m}\neq 0\\
0\quad\quad\text{otherwise}
\end{cases}.
\end{equation}
The entry $R_{jm}$ is nonzero, i.e. $R_{jm}=r_{jm}$, if and only if the metabolite $m$ is an input of the reaction $j$. The algebraic structure of $G$ is thus completely characterized by the network structure, only. \\

The following definition, originally introduced by Brehm and Fiedler \cite{BF18}, presents a central tool for this paper.

\begin{defn}[Child Selections, mothers, children] \label{CSDEF}
A \emph{Child Selection} is an injective map $\mathbf{J}: \textbf{M} \longrightarrow \textbf{E}$, which associates to every metabolite $m \in \textbf{M}$ a reaction $j \in \textbf{E}$ such that $m$ is an input metabolite of reaction $j$.\\
We call the reaction $j=\mathbf{J}(m)$ {a} \emph{child} of $m$, and the metabolite $m=\mathbf{J}^{-1}(j)$ {a} \emph{mother} of the reaction $j$.
\end{defn}


\begin{rmk}
{
It is possible that a metabolite $m$ is an input of $j$ but not a mother of $j$, due to injectivity of Child Selections. Indeed, consider the following example:
}
\begin{equation}
\end{equation}

\vspace{-0.8cm}

\begin{center} 
\includegraphics[scale=0.25]{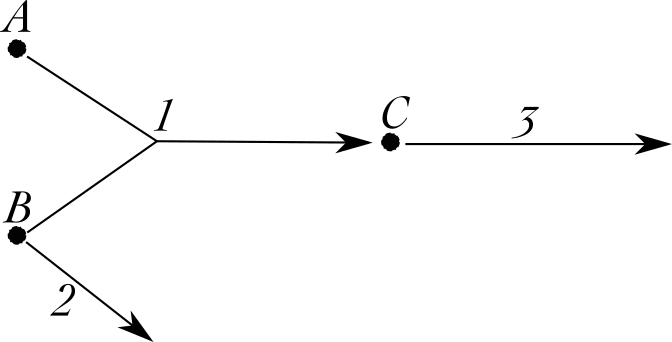}
\end{center}
{
In this minimal case, metabolite $B$ is an input of reaction 1, but never a mother: 1 is the only outgoing reaction from metabolite $A$. Therefore, due to injectivity, $\mathbf{J}^{-1}(1)=A$.}
\end{rmk}

The following Jacobian analysis, based on the Cauchy-Binet formula, is developed from a previous result in \cite{BF18}. 

\begin{prop} \label{det}
Let $G$ be a network Jacobian matrix, in the above settings. Then:
\begin{equation}
\operatorname{det}G=\sum_\mathbf{J} \operatorname{det}S^\mathbf{J} \cdot \prod_{m \in \mathbf{M}} r_{\mathbf{J}(m)m},
\end{equation}
where $S^\mathbf{J}$ is the {$M\times M$} matrix whose $m^{th}$ column is the $\mathbf{J}(m)^{th}$ column of $S$.
\end{prop}

\begin{rmk} \label{dg}
 {Because of the exclusion of explicit autocatalytic reactions from our model, the stoichiometric coefficient $S_{mj}$ of an input metabolite $m$ to reaction $j$ is a priori negative. This implies that any diagonal element of $S^\mathbf{J}$ is always negative: $S^\mathbf{J}_{mm}=S_{m\mathbf{J}(m)}<0$} for any Child Selection $\mathbf{J}$ and any metabolite $m$.
\end{rmk}
\begin{rmk}
If there are no Child Selections at all, then $\operatorname{det}(G)\equiv 0$ for any choice of parameters $r_{jm}$. {For example, this is the case of a network where there are more metabolites than reactions. Another possible example is when two metabolites $m_1$ and $m_2$ are both inputs to one reaction $j$, and they are inputs to no other reaction. In such a situation, due to the injectivity requirement in Definition \ref{CSDEF}, we automatically have no Child Selections and $\operatorname{det}(G)\equiv 0$ for all parameters. Note that, in metabolic networks, this is practically never the case due to the omnipresence of outflow reactions.}
\end{rmk}

Via Proposition \ref{det}, the Jacobian determinant of a metabolic network is a homogen{e}ous multilinear polynomial in the variables $r_{jm}$. The possible sign of the polynomial depends on the signs of the monomial coefficients $\operatorname{det}S^\mathbf{J}$. Hence, it is natural to state the following classification of Child Selections, according to the sign of the  determinant of the reshuffled minor $S^\mathbf{J}$.

\begin{defn}[Child Selection behavior] \label{CSclass}
Let $\mathbf{J}$ be a Child Selection.\\
We say that $\mathbf{J}$ is \emph{bad}, or $\mathbf{J}$  $\emph{ill-behaves}$, if $\operatorname{sign}(\operatorname{det}S^\mathbf{J})=(-1)^{M-1}$.\\
We say that $\mathbf{J}$ is \emph{good}, or $\mathbf{J}$ \emph{well-behaves}, if $\operatorname{sign}(\operatorname{det}S^\mathbf{J}) = (-1)^M$.\\
If $\operatorname{det}S^\mathbf{J}=0$, we say that $\mathbf{J}$ $\emph{zero-behaves}$.
\end{defn}

{The choice of the terminology is consistent. In fact, in a metabolic network context, important classes of Child Selections well-behave, for example all acyclic Child Selections (see Section \ref{IR1}). Moreover, a `stable' Child Selection $\mathbf{J}$, for which all the eigenvalues of $S^\mathbf{J}$ have negative real part, well-behaves. In this sense, a ill-behaving Child Selection is an indication of possible instability.}

\vspace{1cm}
\textbf{Example 1: bad child selection}

This network represents a bad Child Selection:

\begin{equation}\label{usualbad}
\end{equation}

\vspace{-1.2cm}

\begin{center} 
\includegraphics[scale=0.28]{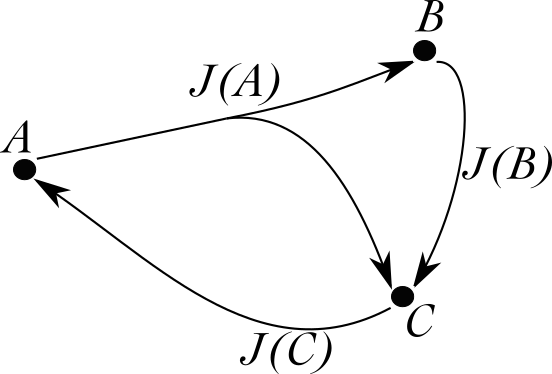}
\end{center}
The stoichiometric matrix associated is
\begin{equation}
S^\mathbf{J}=
\begin{blockarray}{cccc}
 & \mathbf{J}(A) & \mathbf{J}(B) & \mathbf{J}(C)\\
\begin{block}{c[ccc]}
  A & -1 & 0 & 1 & \\
  B & 1 & -1 & 0 &  \\
  C & 1 & 1 & -1& \\
\end{block}
\end{blockarray}\; ,   \;\;
\operatorname{det}S^\mathbf{J}= +1,
\end{equation}
with $\operatorname{sign}(\operatorname{det}S^\mathbf{J})=(-1)^{M-1}=(-1)^{3-1}=+1$.
\vspace{1cm}

\textbf{Example 2: good Child Selection}

This network represents a good Child Selection:

\begin{equation}\label{usualgood}
\end{equation}

\vspace{-1.2cm}

\begin{center} 
\includegraphics[scale=0.28]{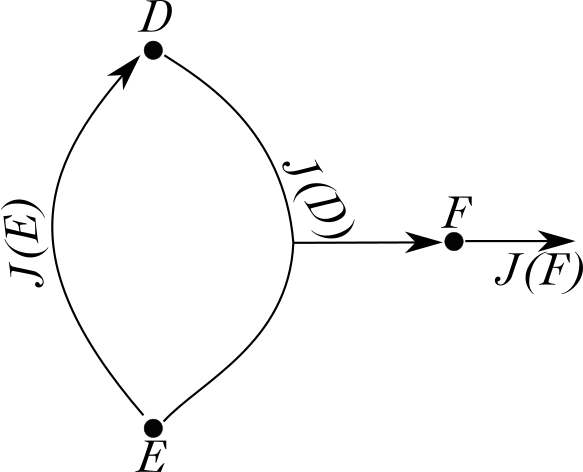}
\end{center}
The stoichiometric matrix associated is
\begin{equation}
S^\mathbf{J}=
\begin{blockarray}{cccc}
 & \mathbf{J}(D) & \mathbf{J}(E) & \mathbf{J}(F)\\
\begin{block}{c[ccc]}
  D & -1 & 1 & 0 \\
  E & -1 & -1 & 0 \\
  F & 1 & 0 & -1 \\
\end{block}
\end{blockarray}\;,   \;\;
\operatorname{det}S^\mathbf{J}= -2,
\end{equation}
with $\operatorname{sign}(\operatorname{det}S^\mathbf{J})=(-1)^M=(-1)^3=-1$.
\vspace{1cm}

At this point, the reader may wonder whether Definition \ref{CSclass} is well-posed in the first place; i.e., whether the behavior of a Child Selection depends on the specific labeling of the network. Section \ref{IR1}, Remark \ref{wprmk}, clarifies this point, assuring the well-posedness of the definition.\\


\section{Main result} \label{MT1}

In metabolic networks, stoichiometric coefficients are mostly 0 and 1. For this reason, we continue here assuming that $S$ has entries $S_{mj}\in \{-1,0,1\}$. By \ref{CBA}, Remark \ref{dg}, then, the diagonal entries of the reshuffled minor $S^\mathbf{J}$ are $S^\mathbf{J}_{mm} \equiv -1$, for any $m$. We refer to the Phd thesis \cite{Vas20} for a more general version of Theorem \ref{MainThm1}, accounting for real stoichiometric coefficients $S_{mj}\in \mathbb{R}$.\\
Here, via a structural analysis of $\operatorname{det}S^\mathbf{J}$, we characterize whether the given Child Selection $\mathbf{J}$ well-behaves or ill-behaves. Note, however, that the importance of the result is mainly revealed in its interpretation, see Section \ref{IR1}.\\

The Leibniz expansion formula for the determinant, applied to $S^\mathbf{J}$, reads
\begin{equation} \label{leibniz}
\operatorname{det}S^\mathbf{J}=\sum_\pi \operatorname{sgn}(\pi) \prod_{m \in \mathbf{M}} S^\mathbf{J}_{\pi(m)m}.
\end{equation}
Again, $\pi$ indicates permutations of $M$ elements and $\operatorname{sgn}(\pi)$ is the signature of $\pi$. Let
\begin{equation}
\mathpzc{E}(\pi):=\operatorname{sgn}(\pi) \prod_{m \in \mathbf{M}} S^\mathbf{J}_{\pi(m)m}
\end{equation}
denote the summand associated to the permutation $\pi$ in the Leibniz expansion. For example, denoting as $\operatorname{\mathpzc{Id}}$ the identity permutation,
\begin{equation}
\mathpzc{E}(\operatorname{\mathpzc{Id}})=(-1)^M.
\end{equation} 
Let $\pi \neq \operatorname{\mathpzc{Id}}$ be a permutation such that $\mathpzc{E}(\pi)\neq 0$. Combinatorially, the permutation $\pi$ can be expressed as the product of $\vartheta$ disjoint permutation cycles $c_i$ of length $l_i>1$, 
\begin{equation}
\pi = \prod_{i=1}^\vartheta c_i.
\end{equation}

\begin{defn}[good/bad-completions, good/bad-cycles] \label{completion}
Given a Child Selection $\mathbf{J}$, we call $\pi$ a \emph{good-completion} if 
\begin{equation}
 \prod_{m:\,\pi(m)\neq m} S^\mathbf{J}_{\pi(m)m} = (-1)^\vartheta.
\end{equation}
We call $\pi$ a \emph{bad-completion} if 
\begin{equation}
 \prod_{m:\,\pi(m)\neq m}S^\mathbf{J}_{\pi(m)m} = (-1)^{\vartheta-1}.
\end{equation}
Again, $\vartheta$ is the number of cycles in the permutation expansion. If $\pi$ consists of a single cycle $c$, i.e. for $\vartheta=1$, we call the good (resp. bad)-completion a \emph{good(resp. bad)-cycle}. 
\end{defn}

We clarify in the next Section \ref{IR1} what a \emph{completion} completes, as it requires some further arguments. Firstly, given the above definition, we state the main result of this section.\\

\begin{thm} \label{MainThm1}
Let $\mathbf{J}$ be a Child Selection and let $\mathbf{\mathpzc{G}}$ and $\mathbf{\mathpzc{B}}$ be the number of good and bad completions, respectively. Then, in the sense of Definition \ref{CSclass},
\begin{enumerate}[itemsep=0pt]
\item The Child Selection $\mathbf{J}$ well-behaves if $\mathbf{\mathpzc{G}} > \mathbf{\mathpzc{B}}-1$.
\item The Child Selection $\mathbf{J}$ ill-behaves if $\mathbf{\mathpzc{G}}<\mathbf{\mathpzc{B}}-1$.
\item The Child Selection $\mathbf{J}$ zero-behaves if $\mathbf{\mathpzc{G}} = \mathbf{\mathpzc{B}}-1$.
\end{enumerate}
\end{thm}

For a Child Selection $\mathbf{J}$, Theorem \ref{MainThm1} characterizes the sign of $\operatorname{det}S^\mathbf{J}$ in terms of permutation cycles of the Leibniz expansion \eqref{leibniz}. The following Section \ref{IR1} relates these permutation cycles to certain cycles in the network.

\section{Graphical interpretation of the result} \label{IR1}

\begin{figure}
\begin{center}

\begin{TAB}(r,2cm,2cm)[5pt]{|c|c|c|c|c|}{|c|c|c|}
Examples & 

\underline{\textit{BIOLOGICAL}} & 

\begin{minipage}{.2\textwidth}
\underline{\textit{MR-GRAPH}}\\ 
\textit{Biological} \\
   \end{minipage}&
  
\begin{minipage}{.2\textwidth}
\underline{\textit{MR-GRAPH}}\\ 
\textit{Combinatorial} \\
   \end{minipage}
   &

   \underline{\textit{MATRIX}}
   
   \\

1 &
 \begin{minipage}{.2\textwidth}
      \includegraphics[scale=0.21]{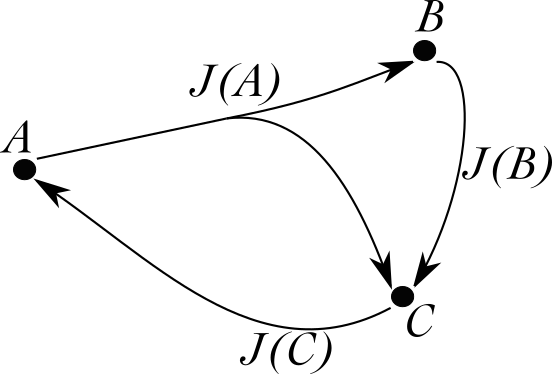}
  \end{minipage}
  
  & 
   \begin{minipage}{.2\textwidth}
      \includegraphics[scale=0.21]{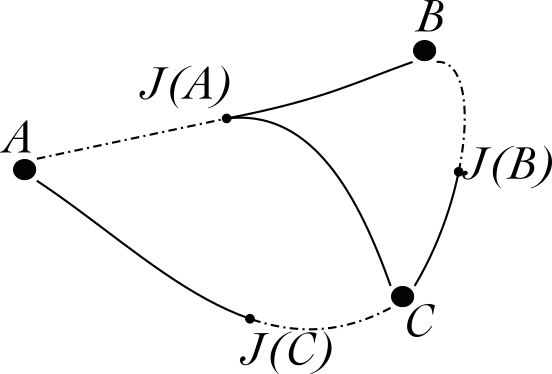}
  \end{minipage}
  
  &
   \begin{minipage}{.2\textwidth}
      \includegraphics[scale=0.25]{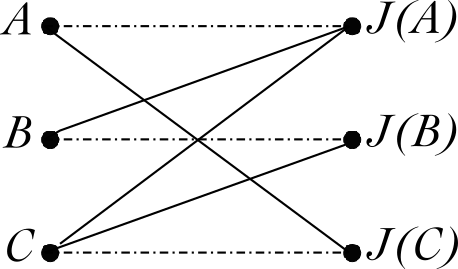}
  \end{minipage}
  
  &
 
\begin{blockarray}{cccc}
 & $\mathbf{J}(A)$ & $\mathbf{J}$(B) & $\mathbf{J}$(C)\\
\begin{block}{c[ccc]}
  A & -1 & 0 & 1 \\
  B & 1 & -1 & 0 \\
  C & 1 & 1 & -1\\
\end{block}
\end{blockarray}\\
 2 & 
\begin{minipage}{.2\textwidth}
      \includegraphics[scale=0.2]{ex2bio.png}
  \end{minipage} & 
\begin{minipage}{.2\textwidth}
      \includegraphics[scale=0.2]{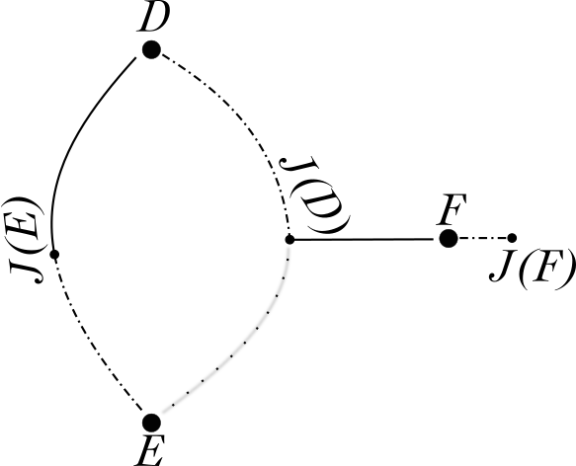}
  \end{minipage} &
  \begin{minipage}{.2\textwidth}
      \includegraphics[scale=0.25]{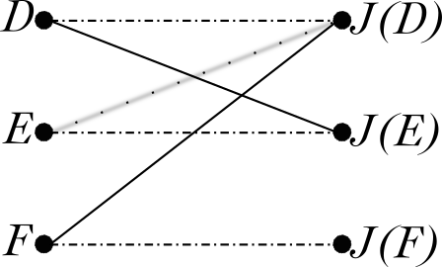}
  \end{minipage}
 &
\begin{blockarray}{cccc}
 & $\mathbf{J}$(D) & $\mathbf{J}$(E) & $\mathbf{J}$(F)\\
\begin{block}{c[ccc]}
  D & -1 & 1 & 0 \\
  E & -1 & -1 & 0 \\
  F & 1 & 0 & -1 \\
\end{block}
\end{blockarray}

   \\   
\end{TAB}
\caption{For the two examples of Child Selections, four different ways of representation: biological, $MR$-Graph (in a biological shape), $MR$-graph (in a combinatorial shape), matrix. Note that, when labeled, the four representations are equivalent. In the $MR$-graphs,  \emph{negative edges \textbf{J}-selected} are indicated with a dotted-dashed line, the sparse dotted {grey} line indicates \emph{negative edges not \textbf{J}-selected}, the continuous black line indicates \emph{positive edges}. In the combinatorial shape, the edges \textbf{J}-selected are the horizontal ones. Example 1 possesses two completion cycles: $c_1=A-\mathbf{J}(A)-C-\mathbf{J}(C)-A$ and $c_2=A-\mathbf{J}(A)-B-\mathbf{J}(B)-C-\mathbf{J}(C)-A$, both bad. Example 2 possesses only one good completion cycle: $c=D-\mathbf{J}(D)-E-\mathbf{J}(E)-D$.  Consequently, Example 1 represents a bad Child Selection, and Example 2 represents a good one.} 
\label{ESEMPI}
\end{center}
\end{figure}

The \emph{Metabolite-Reaction graph} ($MR$-graph) is an undirected bipartite graph with two sets of vertices given by the metabolites $m_1,...,m_M$ and the reactions $j_1,...,j_E$, respectively. For a metabolite $m$ participating in a reaction $j$, edges $e=(m,j)$ are adjacent to a metabolite vertex $m$ and a reaction vertex $j$. With this construction, edges in the $MR$-graph are in a one-to-one relation with the nonzero entries of stoichiometric matrix $S$. In particular, with $S_{mj}<0$ of 
\begin{equation}
j: \quad m+... \underset{j}{\longrightarrow} ...
\end{equation} 
in mind, we call \emph{negative} the edges $e=(m,j)$ where $m$ is input to $j$. Conversely, with $S_{mj}>0$ of 
\begin{equation}
j: \quad ... \underset{j}{\longrightarrow} m+ ...
\end{equation}
in mind, we call \emph{positive} the edges $e=(m,j)$ where $m$ is output to $j$. See figure \ref{ESEMPI}  for a comparison between different kinds of representation graphs for the same network. Under the name Species-Reaction graph ($SR$-graph), this was considered by \cite{CraFei06} and others.\\

We proceed with two definitions and a proposition.

\begin{defn}[\textbf{J}-selected edges]
For any Child Selection $\mathbf{J}$, we call the negative edges $e=(m, \mathbf{J}(m))$ in the $MR$-graph $\mathbf{J}$\emph{-selected}. 
\end{defn}

\begin{rmk} 
In particular, $\mathbf{J}$-selected edges are such that the corresponding stoichiometric entry lies on the diagonal of $S^\mathbf{J}$. 
\end{rmk}
\begin{rmk} \label{injadj}
Injectivity of a Child Selection $\mathbf{J}$ directly implies that two $\mathbf{J}$-selected edges $e_1$ and $e_2$ never share a vertex in the $MR$-graph. 
\end{rmk}

\begin{defn} [Completion Cycle]
For a Child Selection $\mathbf{J}$, a completion cycle in the $MR$-graph is a cycle of length $2l$, $\ell \le M$, such that $\ell$ edges are $\mathbf{J}$-selected.
\end{defn}
\begin{rmk}
Equivalently, because of Remark \ref{injadj} above, a completion cycle is a cycle in the $MR$-graph of length $2l$, such that $\ell$ $\mathbf{J}$-selected edges alternate with $\ell$ non $\mathbf{J}$-selected edges.
\end{rmk}

\begin{prop} \label{crucialinterpret}
For any given Child Selection $\mathbf{J}$, there is a one-to-one correspondence between completion cycles and nonzero permutation cycles; i.e., cycles $c$ such that 
\begin{equation}
\prod_{m:\,c(m)\neq m} S^\mathbf{J}_{c(m)m}\neq 0.
\end{equation}
\end{prop}

It is now clear what the word `completion' refers to. In fact, any completion cycle is constructed by completing $\ell$ $\mathbf{J}$-selected elements to a cycle of length $2l$ in the $MR$-graph. In this sense, a good-completion $\pi = \prod_{i=1}^\vartheta c_i$ can be seen, in the $MR$-graph, as a collection of $\vartheta$ non-intersecting completion cycles $c_i$, such that the number of good-cycles has the same parity of $\vartheta$. Respectively, a bad-completion is a collection of $\vartheta$ non-intersecting completion cycles, such that the number of good-cycles has opposite parity of $\vartheta$.

{
\begin{rmk}
From Definition \ref{completion} and Proposition \ref{crucialinterpret} it follows that, for a completion cycle in the $MR$-graph, the parity of \emph{the negative edges not \textbf{J}-selected} characterizes the cycle being good or bad. In fact: an odd number of negative edges not \textbf{J}-selected corresponds to a good completion cycle, an even number corresponds to a bad one. 
\end{rmk}}

\begin{rmk} \label{wprmk}
Obviously, being a network structure, the definition of a completion cycle does not depend on the specific labeling of the network. Proposition \ref{crucialinterpret}, together with Theorem \ref{MainThm1}, guarantees in particular that also Definition \ref{CSclass} does not depend on any labeling and, thus, is well-posed.
\end{rmk}

Finally, we list some consequences of Theorem \ref{MainThm1}, useful for applications. 

\begin{cor}[Examples of application]\label{exofapp}
The following statements hold true:
\begin{enumerate}[itemsep=0pt]
\item Acyclic Child Selections well-behave;
\item A Child Selection possessing one good-cycle, and no other completion cycles, well-behaves;
\item A Child Selection possessing one bad-cycle, and no other completion cycles, zero-behaves;
\item A Child Selection possessing two intersecting bad-cycles, and no other completion cycles, ill-behaves;
\item Any nonzero Child Selection of a network which possesses only monomolecular reactions and one single bimolecular reaction
\begin{equation} \label{stupeq}
\tilde{j}:\quad A+B \underset{\tilde{j}}{\longrightarrow} C
\end{equation}
well-behaves.
\end{enumerate}
\end{cor}

In particular, the bad Example 1 falls into the casuistry of the point (4) of Corollary \ref{exofapp}, and the good Example 2 into the casuistry of point (2) of Corollary \ref{exofapp}. Any Child Selection with a single monomolecular cycle is an example for (3). Such list can be of course enlarged for a given network of application, to easily detect good and bad Child Selections for bifurcation analysis, as discussed in the continuation of this paper.\\

\section{Hunting saddle-node bifurcations} \label{SNB}

Here, we give a simple network condition under which there is the possibility, for certain parameters, of a saddle-node bifurcation of equilibria. We also identify bifurcation parameters responsible for the change of sign of the determinant and consequent change of stability of any equilibrium.

\begin{thm}[Change of Stability] \label{CoS}
Suppose there exist two Child Selections $\mathbf{J}_1$, $\mathbf{J}_2$, and a metabolite $m_b$, such that $\mathbf{J}_1(m_b)\neq \mathbf{J}_2(m_b)$ and $\mathbf{J}_1(m) = \mathbf{J}_2(m)$ for any $m\neq m_b$. Assume moreover that $\mathbf{J}_1$ well-behaves and $\mathbf{J}_2$ ill-behaves. Then the Jacobian determinant of $G$ takes the form
\begin{equation}
\operatorname{det}G=(\mathpzc{a}  r_{\mathbf{J}_1(m_b)m_b}-\mathpzc{b} r_{\mathbf{J}_2(m_b)m_b}) r_{\mathbf{J}_1(m_{n+1})m_{n+1}} ... r_{\mathbf{J}_1(m_m)m_m}\quad+\quad ... \quad ,
\end{equation}
where the omitted terms can be chosen arbitrarily small, and $\mathpzc{a},\, \mathpzc{b}$ are coefficients of the same sign. In particular, the parameter
\begin{equation}\label{csi}
\xi=\mathpzc{a} r_{\mathbf{J}_1(m_b)m_b}-\mathpzc{b}  r_{\mathbf{J}_2(m_b)m_b}
\end{equation}
may serve as a bifurcation parameter for the bifurcation of nontrivial equilibrium solutions of the system \eqref{eq}.
\end{thm}

Theorem \ref{CoS} does not require fixing an equilibrium, and it holds in more general settings. Nevertheless, for simplicity, we are thinking here of an equilibrium situation. The parameter $\xi=\mathpzc{a}  r_{\mathbf{J}_1(m_b)m_b}-\mathpzc{b}  r_{\mathbf{J}_2(m_b)m_b}$ is `localized' in a single metabolite $m_b$. In fact, the change of stability is driven by the difference between the derivatives with respect to the same metabolite $m_b$ of the reaction rates of two child reactions of $m_b$ itself. This suggests a simple biological scheme for controlling the unstable dimension of the equilibrium.\\

\section{A case study: the central metabolism of {\textit{E. coli}}} \label{TCAC}

The central metabolism of {\textit{E. coli}} (figure \ref{anna}) consists of different and interconnected parts. In particular, the upper part comprises the so-called \emph{Pentose phosphate pathway} and the \emph{Glycolysis}. The bottom `cyclic' part includes the \emph{Tricarboxylic acid cycle} and the \emph{Glyoxylate cycle}. We skip here more detailed biological explanation \cite{Al03, Lod08}.\\

In this section, we analyze the network of the central metabolism of {\textit{E. coli}} of figure \ref{anna}. This network representation is primarily based on the original model proposed by Ishii et al. in \cite{Ish07} with the modifications suggested by Nakahigashi et al. in \cite{Nak+09}. Note that, in biology papers, `obvious' outflow exit reactions, such as d1 - d6 shown here, are frequently omitted. For our mathematical analysis, however, we are bound to include them. These reactions are the only outgoing reactions of their input metabolites. Their omission would result in an infinite production of their input metabolites and in a mathematical degeneracy of the network.\\

The network possesses 30 metabolites and 58 reactions. The number of Child Selections is on the order of $10^7$. Nevertheless, we can provide interesting biological insights without computing such a huge amount of Child Selections. In the same spirit as Section \ref{SNB}, and along its lines, we find two Child Selections $\mathbf{J}_1$ and $\mathbf{J}_2$ with opposite behavior, such that $\mathbf{J}_1(m_b) \neq \mathbf{J}_2(m_b)$ for one single metabolite $m_b$ and $\mathbf{J}_1(m) = \mathbf{J}_2(m)$ for all other metabolites $m \neq m_b$. This situation, via Theorem \ref{CoS}, provides a bifurcation parameter responsible for a change of sign in the Jacobian determinant and possible consequent saddle-node bifurcations of equilibria.\\

To find the two Child Selections $\mathbf{J}_1$ and $\mathbf{J}_2$ as above, we start by imposing certain child reactions $j$ to certain mother metabolites $m$. We do this arbitrarily, and only for sake of exemplification. Many other choices and analogous constructions are, of course, possible.\\ 
Let us {consider only Child Selections associating the metabolites $PEP$, $PYR$, and $CO2$ to their respective exit reactions, that is:} 
\begin{enumerate}[itemsep=0pt]
\item $\mathbf{J}_1(PEP)=\mathbf{J}_2(PEP)=dd6$;
\item $\mathbf{J}_1(PYR)=\mathbf{J}_2(PYR)=dd7$;
\item $\mathbf{J}_1(CO2)=\mathbf{J}_2(CO2)=d6$.
\end{enumerate}

{For any Child Selection satisfying these constraints, we can} consider the upper part (Pent. Phosph. Pathway - Glycolysis) and the bottom bart (Tricarboxylic acid cycle - Glyoxylate cycle) as separate and independent. In fact, any Child Selection $\mathbf{J}$ satisfying the above constraints 1-3, identifies reshuffled minors $S^{\mathbf{J}}$, which are block diagonal. This shows that certain qualitative arguments on the dynamics of the central metabolism may be inferred, separately, from the biological components of the network. For example, for a block diagonal Jacobian matrix in our settings, indeterminate sign determinant of one block trivially implies indeterminate sign determinant for the entire matrix. In particular, we may concentrate on the bottom part of the network, only assuming that $\mathbf{J}_1 = \mathbf{J}_2$ in the upper part.\\

\begin{figure}
\begin{center}
\begin{tabular}{ c c } 
\underline{\textit{BIOLOGICAL}} &
\underline{\textit{MR-GRAPH}} (\textit{Biological}) \\

     \includegraphics[scale=0.9]{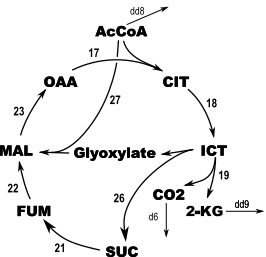}&

     \includegraphics[scale=0.9]{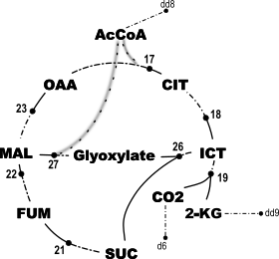}\\
\end{tabular}
\caption{{This subnetwork focuses on the Tricarboxylic cycle and the Glyoxylate cycle. It possesses only two Child Selections, depending on the child of ICT. In fact, note that reaction 17 is the single child of OAA and reaction 27 is the single child of Glyoxylate. Due to injectivity of Child Selections, then, $AcCoA$ also has only one child: reaction $dd8$.}}
 \label{twonet}
\end{center}
\end{figure}

In figure \ref{twonet}, we have depicted the chosen subnetwork possessing precisely only the two chosen Child Selections $\mathbf{J}_1$ and $\mathbf{J}_2$. We identify the metabolite $m_b$, such that $\mathbf{J}_1 (m_b) \neq \mathbf{J}_2 (m_b)$, as $ICT=m_b$. Indeed, any other metabolite $m \neq ICT$ possesses a single child $\mathbf{J}_1 (m) = \mathbf{J}_2 (m)$, and only $ICT$ possesses two child reactions: reaction 19 and reaction 26. Let us call $\mathbf{J}_{1}$ the Child Selection such that $\mathbf{J}_{1}(ICT)=19$ and $\mathbf{J}_{2}$ the Child Selection such that $\mathbf{J}_{2}(ICT)=26$. With this choice of Child Selections, metabolites \emph{Lactate, Acetate}, and \emph{Ethanol} result disconnected from the rest of the network and have consequently been omitted here.\\

By looking at the $MR$-graph representation, we can easily conclude that $\mathbf{J}_1$ well-behaves and $\mathbf{J}_2$ ill-behaves.\\ 
Note indeed that $\mathbf{J}_{1}$ does not contain any completion cycle, and therefore well-behaves. In fact, this Child Selection contains only one network cycle $c=MAL-23-OAA-17-AcCoa-27-MAL$, which is not a completion cycle as the edge $AcCoa-27$ is not $\mathbf{J}_{1}$-selected.\\
On the other hand, the completion cycles structure of $\mathbf{J}_{2}$ is identical to Example 1, which had provided a simple and recognizable pattern of an bad Child Selection. In fact, this Child Selection possesses only two bad completion cycle $c_1$ and $c_2$:
\begin{enumerate}[itemsep=0pt]
\item $c_1=ICT-26-Glyoxylate-27-MAL-23-OAA-17-CIT-18-ICT$;
\item $c_2=ICT-26-SUC-21-FUM-22-MAL-23-OAA-17-CIT-18-ICT$.
\end{enumerate}
Since it possesses only two intersecting bad completion cycles, the Child Selection $\mathbf{J}_{26}$ ill-behaves.\\
In particular, in accordance to Theorem \ref{CoS}, the parameter
\begin{equation}
\xi=r_{19m_b} - r_{26m_b},\text{ (where }m_b=ICT\text{)},
\end{equation}
controls a change of sign of the Jacobian determinant of the entire system, for a certain region of parameters.

\begin{rmk} 
The choice of reaction $19$ and $26$ basically highlights the difference between the Tricarboxylic acid cycle (reaction 19) and the Glyoxylate cycle (reaction 26). Our analysis suggests how the control of certain dynamical properties of the metabolism of a cell is related to its network structure.

\end{rmk}

\section{Discussion}\label{discussion}

In this paper, we have presented a new approach to address questions related to the Jacobian determinant for dynamical systems arising from metabolic networks. The idea of our approach relies on considering the partial derivatives $r_{jm}(x^*)$ independent, at a fixed equilibrium $x^*$, from the value of the equilibrium $x^*$ itself and from the value $r_j(x^*)$ attained by the reaction rates at the equilibrium. Importantly, this allows us to separate the question of the existence of the equilibrium, which concern the reaction rates $r_j$, from the questions related to the stability of the equilibrium, which concern the derivatives $r_{jm}$.\\ 
In particular, we interpret the Jacobian determinant as a homogen{e}ous multilinear polynomial in the variables $r_{jm}$. The coefficients of each monomial are determinants of reshuffled stoichiometric minors $S^\mathbf{J}$, and we have given network conditions to establish the sign of $\operatorname{det}S^\mathbf{J}$.\\

{Mass action kinetics forbids considering the partial derivatives $r_{jm}$ as free parameters, at an equilibrium. Still, the present structural analysis of the Jacobian determinant holds identically also for the mass action case. In particular, a Jacobian determinant of fixed sign is still a sufficient condition to exclude saddle-node bifurcations also under mass action kinetics, and the copresence of both good and bad Child Selections is therefore still a necessary condition. However, in the case where the sign of the Jacobian depends on the parameters $r_{jm}$, we are not able to directly conclude on equilibria bifurcations, for mass action systems, given that the existence of the equilibrium depends as well on related parameters. Hence, the consistency between the parameter constraints coming from the existence of the equilibrium and the constraints coming from the Jacobian determinant must be additionally checked. This introduces a further difficulty in bifurcation analysis under mass action kinetics.\\}

Most of the arguments can be lifted to any dynamical system of the form
$$\dot{x}=Ar(x),$$
where $A$ is any real matrix (see \cite{Vas20}). However, the efficacy of this approach is most exploited for the case where $A=S$ is a sparse matrix with low integer entries, as it is the case of stoichiometric matrices of metabolic networks. In this context, our approach provides bifurcation patterns and biological intuition.\\
Theorem \ref{CoS} prescribes how to find parameters for a change of stability of equilibria. In particular, we shall find one good Child Selection $\mathbf{J}_1$ and one bad Child Selection $\mathbf{J}_2$ such that $\mathbf{J}_1(m_b)\neq \mathbf{J}_2(m_b)$ for a single metabolite $m_b$, and $\mathbf{J}_1(m)=\mathbf{J}_2(m)$ for all other metabolites $m\neq m_b$.\\
We point out that finding good Child Selections is reasonably easy, in metabolic networks. As an example, the sparsity of the stoichiometric matrix $S$ as well as the overall presence of outflows make the task of finding acyclic Child Selections simple and, by Corollary $\ref{exofapp}$, acyclic Child Selections are always good. An analogous argument supporting the predominance of good Child Selections follows by the fact that many reactions in a metabolic networks are monomolecular, always via $\ref{exofapp}$.\\
Morally, then, the fewer bad Child Selections are the important ones to be found in the quest for parameter areas where bifurcations and changes of stability take place. This last observation underlines the importance to find simple and recognizable patterns of bad Child Selections. Example 1 of this paper is possibly one of those examples. Designed as the simplest example of a bad Child Selection on three metabolites, it turned out to possess a cycle structure that was easy to recognize in the real example of the central carbon metabolism of {\textit{Escherichia coli}}.\\

{Previous studies addressed equilibria bifurcations in the central carbon metabolism of \textit{E. coli}.\\
In \cite{VL06} the authors have numerically studied bifurcations in the restricted \textit{E. coli} model proposed by Chassagnole and coauthors \cite{Chass2002}. They have focused on saddle-node bifurcations and Hopf bifurcations. In a similar flavor to the present paper, their bifurcation parameters are associated to individual reactions $j$. In particular, the rate $r_j$ is written as $r_j=k_j g(x)$, where $k_j$ is a scalar and $g(x)$ the given kinetics. They use the scalars $k_j$ as bifurcation parameters. However, changing this parameter also influences the reference equilibrium $x^*$. For mathematical simplicity, in contrast, our approach fixes $x^*$ and changes the partial derivative $r_{jm}$, only. With these adaptations in mind, our analytical results match their numerical explorations well. For example, consider the following subnetwork from their restricted model:}
\begin{equation}\label{other}
\end{equation}

\vspace{-0.8cm}

\begin{center} 
\includegraphics[scale=0.25]{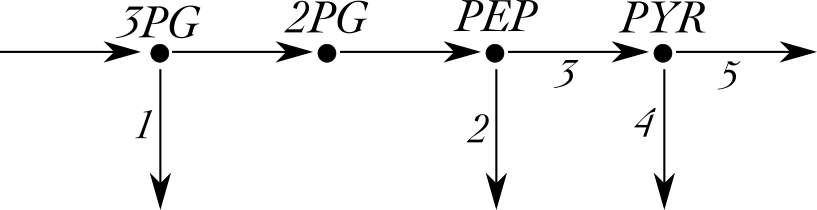}
\end{center}

{Here reactions $1,2,4,5$ are outflow reactions, and $PYR$ is not an input to any other reaction. Their bifurcation diagram for the two parameters $(k_1,k_2)$ is similar, qualitatively, as for $(k_1,k_3)$. Qualitatively similar diagrams also appear when considering the couples $(k_1,k_4)$ and $(k_1,k_5)$. This is consistent with our Child Selections point of view. Indeed: assume that $\mathbf{J}$ is any good Child Selection with $\mathbf{J}(PEP)=2$. Then, it is good also the Child Selection $\tilde{\mathbf{J}}$ such that $\tilde{\mathbf{J}}(PEP)=3$ and $\tilde{\mathbf{J}}(m)=\mathbf{J}(m)$ for any other metabolite $m$, because the cycle structure remains identical: the subnetwork \eqref{other} is acyclic. The same argument holds also for the other two `sisters', reactions 4 and 5. Our analysis consistently suggests that the two parameters $k_2$ and $k_3$ ($k_4$ and $k_5$, respectively) are bifurcation parameters that can be used in an analogous way, qualitatively.\\
The absence of equilibria bifurcations in the \textit{E. coli} network, both for mass action and Michaelis-Menten kinetics, has been claimed in \cite{OTM18}. Curiously, that study misses the sign changes of the Jacobian determinant, which we encounter for the Michaelis-Menten case.\\
Mathematically, global Hopf bifurcation of time periodic oscillations of the \textit{E. coli} network, and specifically for the citric acid cycle part, has also been established along the lines of our present analysis \cite{F19}.\\}

{There are several applications of the theory presented in this paper. Firstly, our theory suggests promising \emph{bifurcations} parameters based on network information, only. In large systems, finding bifurcation parameters for numerical simulations is not a simple task. Our recipe of finding two Child Selections, one good and one bad, serves as an aid for this purpose. However, for realistic explanation of metabolic phenomena the mathematical approach should best be combined with deep biological insight. Further, the tools developed here also apply to \emph{sensitivity} analysis of equilibria. See the thesis \cite{Vas20} where - in the same setting - responses to perturbations of reaction rates and metabolites concentration have been studied. The responses may have a definite sign for all parameters, or that sign may depend on the parameters choice, pretty much in the same spirit as presented in the present paper. The case in which the sign of the responses depends on the parameters suggests the possibility of \emph{controlling} the sign, as the responses may be positive or negative depending on the parameters. For example, note that stable equilibria $x^*$ require the sign of the Jacobian to be $(-1)^M$ for $M$ metabolites. Such control may therefore indicate external switches between different stable metabolic pathways, or may even induce stem cell differentiation.}

\section{Proofs}\label{proofs}

This Section is devoted to the proofs of our results. We start with Proposition \ref{det}.

\proof[Proof of Proposition \ref{det}]
We apply the Cauchy-Binet formula on $G=SR$ to obtain:
\begin{equation}
\begin{split}
\operatorname{det} G =& \sum_{|\mathcal{E}|=M} \operatorname{det}S^\mathcal{E} \cdot \operatorname{det}R_\mathcal{E} = \sum_{|\mathcal{E}|=M} \operatorname{det}S^\mathcal{E} \; (\sum_{\pi} \operatorname{sgn}(\pi) \cdot \prod_{m\in \mathbf{M}}r_{\pi(m)m}).
\end{split}
\end{equation}
Here $\pi$ indicates a permutation of $M$ elements and $\operatorname{sgn}(\pi)$ is the signature (or parity) of $\pi$. Note that $\prod_{m\in \mathbf{M}} r_{\pi(m)m}\neq 0$ if and only if there is an associated Child Selection $\mathbf{J}$ such that $r_{\mathbf{J}(m)m}=r_{\pi(m)m}$, for every $m$. In particular, the sum runs non trivially only for the selected minors $S^\mathcal{E}$ such that the set $\mathcal{E}$ is the image of $\mathbf{M}$ through a Child Selection $\mathbf{J}$. Now,
\begin{equation}
\begin{split}
\sum_{|\mathcal{E}|=M} \operatorname{det}S^\mathcal{E} \; (\sum_{\pi} \operatorname{sgn}(\pi) \cdot \prod_{m\in \mathbf{M}} r_{\pi(m)m})=& \sum_{\mathcal{E}=\mathbf{J}(\mathbf{M})}  \operatorname{det}S^\mathcal{E} \; (\sum_\mathbf{J} \operatorname{sgn}(\mathbf{J}) \cdot \prod_{m\in \mathbf{M}} r_{\mathbf{J}(m)m})\\
=&\sum_\mathbf{J} \operatorname{det}S^\mathbf{J} \cdot \prod_{m\in \mathbf{M}} r_{\mathbf{J}(m)m}.
\end{split}
\end{equation}
Last step is the observation:
\begin{equation}
\operatorname{det}S^{\mathcal{E}=\mathbf{J}(\mathbf{M})} \cdot \operatorname{sgn}(\mathbf{J}) = \operatorname{det}S^\mathbf{J}.
\end{equation}
\endproof

We proceed with the proof of main technical result, Theorem \ref{MainThm1}.

\proof[Proof of Theorem \ref{MainThm1}]
The proof follows an idea of Banaji and Craciun \cite{BaCra10}. Firstly note:
\begin{equation}
\begin{split}
(\operatorname{det}S^\mathbf{J}) (-1)^M =& (\operatorname{det}S^\mathbf{J}) \; \mathpzc{E}(\operatorname{\mathpzc{Id}})\\
=&\sum_\pi \mathpzc{E}(\pi)\mathpzc{E}(\operatorname{\mathpzc{Id}})\\
=& 1 + \sum_{\pi \neq \operatorname{\mathpzc{Id}}}\mathpzc{E}(\pi)\mathpzc{E}(\operatorname{\mathpzc{Id}}).
\end{split}
\end{equation}
Let $h$ be the number of elements $m$ such that $\pi(m)\neq m$. That is, $h$ is the number of elements of $\pi$ which are not fixed points of the permutation, but belong to a permutation cycle.
\begin{equation} \label{asdr}
\begin{split}
\mathpzc{E}(\pi)\mathpzc{E}(\operatorname{\mathpzc{Id}})=& \operatorname{sgn}(\pi) \left(\prod_{m \in \mathbf{M}} S^\mathbf{J}_{\pi(m)m}\right) \; \operatorname{sgn}(\operatorname{\mathpzc{Id}}) \prod_{m \in \mathbf{M}} S^\mathbf{J}_{mm}\\
=&\left(\prod_{m:\,\pi(m)=m} (S^\mathbf{J}_{mm})^2\right) \; \prod_{i=1}^\vartheta \operatorname{sgn}(c_i) \; \left(\prod_{m:\,\pi(m)\neq m}(S^\mathbf{J}_{\pi(m)m}S^\mathbf{J}_{mm})\right)\\
=&\, (-1)^h \prod_{i=1}^\vartheta \operatorname{sgn}(c_i) \prod_{m:\,\pi(m)\neq m} S^\mathbf{J}_{\pi(m)m}\\
=&\, (-1)^\vartheta \prod_{m:\,\pi(m)\neq m} S^\mathbf{J}_{\pi(m)m}.
\end{split}
\end{equation}
The steps above are made noting that $(S^\mathbf{J}_{mm})^2\equiv 1$, for any $m$ and that, for a cycle $c$ of length $\ell$, $\operatorname{sgn}(c) (-1)^\ell  = -1$. We conclude the proof by observing that 
\begin{equation}
(-1)^\vartheta \prod_{m:\,\pi(m)\neq m} S^\mathbf{J}_{\pi(m)m}=1\text{ (-1, respectively)} 
\end{equation}
if $\pi$ is a good-completion (bad-completion, respectively).  This yields to the identity
\begin{equation} \label{detsJ}
\operatorname{det}S^\mathbf{J} (-1)^M = 1 + \mathbf{\mathpzc{G}} - \mathbf{\mathpzc{B}},
\end{equation}
which proves the Theorem.
\endproof

The interpretation of the above result has been discussed in the Proposition \ref{crucialinterpret} and consequent Corollary \ref{exofapp}.

\proof[Proof of Proposition \ref{crucialinterpret}]
Let us consider the computation \eqref{asdr} in the proof of Theorem \ref{MainThm1}. Consider, for simplicity, a single cycle permutation $\pi=c$ and concentrate on the expression
\begin{equation} \label{CICLO}
\prod_{m:\,c(m)\neq m} S^\mathbf{J}_{c(m)m}S^\mathbf{J}_{mm}.
\end{equation}
Note that the diagonal elements $S^\mathbf{J}_{mm}$ and $S^\mathbf{J}_{c(m)c(m)}$ represents $\mathbf{J}$-selected edges. $S_{c(m)m}$ shares the same column (i.e., reaction vertex) with $S^\mathbf{J}_{mm}$ and the same row (i.e., metabolite vertex) with $S^\mathbf{J}_{c(m)c(m)}$. Following the order of the cycle $c$ in the Expression \eqref{CICLO} leads to the desired identification.
\endproof 

\proof[Proof of Corollary \ref{exofapp}] \textcolor{white}{easily} \\

\begin{enumerate}
\item Acyclic Child Selections do not possess any completion cycle and consequently do not possess any completion. By Theorem \ref{MainThm1}, $\mathpzc{G}=\mathpzc{B}=0$, and thus acyclic Child Selection{s} well-behave.
\item Analogously, a Child Selection $\mathbf{J}$ possessing only one good cycle possesses only one good completion. By Theorem \ref{MainThm1}, $\mathpzc{G}=1$, $\mathpzc{B}=0$, and thus $\mathbf{J}$ well-behaves.
\item As in (2), with $\mathpzc{G}=0$, $\mathpzc{B}=1$. Via Theorem \ref{MainThm1} we have the conclusion.
\item Since a completion is a collection of non-intersecting completion cycles, then two intersecting completion cycles implies two different completions. As above, both completions must be bad. Hence, again by Theorem \ref{MainThm1},
$\mathpzc{G}=0$, $\mathpzc{B}=2$.
\item Note that a non-zero Child Selection $\mathbf{J}$, $\operatorname{det}S^\mathbf{J} \neq 0$, a priori, cannot contain monomolecular cycles of the form:
$$m_0 \longrightarrow m_1 \longrightarrow ... \longrightarrow m_0,$$
where each reaction is monomolecular, i.e.
$$m_i \longrightarrow m_{i+1}.$$
Therefore, for a non-zero Child Selection $\mathbf{J}$ possessing only monomolecular reactions and a single bimolecular reaction $\tilde{j}$ of the form \eqref{stupeq}, either $\mathbf{J}$ does not contain any completion cycle and consequently well-behaves, due to point (1) of this corollary, or $\mathbf{J}$ contains completion cycles in the $MR$-graph which include $\tilde{j}$ as a reaction vertex. In particular all completion cycles intersects.\\
Without loss of generalities, assume $\mathbf{J}^{-1}(\tilde{j})=A$. We only have two possibilities for a completion cycle: a bad completion cycle containing the two adjacent edges $A-\tilde{j}-C$ or a good completion cycle containing the two adjacent edges $ A-\tilde{j}-B  $. By assumption on the network, the completion cycles run through monomolecular reactions only, except from $\tilde{j}$. Because they do intersect, we have at most one bad completion, $\mathpzc{G}\le 1$, and one good completion, $\mathpzc{B} \le 1$. Hence, $1-\mathpzc{G}+\mathpzc{B} \ge 0$ and by Theorem \ref{MainThm1} the Child Selection $\mathbf{J}$ does not ill-behave. Since, by assumption, $\mathbf{J}$ does not zero-behave, it must well-behave.
\end{enumerate}
\endproof

The last remaining proof is Theorem \ref{CoS}. We introduce some concepts, first.\\ 
The set of Child Selections $\{\mathbf{J}\}$ carries a natural integer-valued distance $d$. 
\begin{defn}
Let $\mathbf{J}_1$, $\mathbf{J}_2$ be two Child Selections. We define the distance $d(\mathbf{J}_1,\mathbf{J}_2)$ as the number of metabolites $m \in \mathbf{M}$ such that $\mathbf{J}_1(m) \neq \mathbf{J}_2(m)$.
\end{defn}

It is straightforward to verify that $d$ is a distance on the set of Child Selections. We consider now Child Selections at distance $d=1$. These are Child Selections $\mathbf{J}_1$, $\mathbf{J}_2$ such that $\mathbf{J}_1(m_b)\neq \mathbf{J}_2(m_b)$ for a single metabolite $m_b$ and $\mathbf{J}_1(m)=\mathbf{J}_2(m)$ for any $m\neq m_b$ different from $m_b$, as in Theorem \ref{CoS}. Clearly:
\begin{equation}
\begin{split}
&\operatorname{det}S^{\mathbf{J}_1}\prod_m r_{{\mathbf{J}_1}(m)m}+\operatorname{det}S^{\mathbf{J}_2}\prod_m r_{{\mathbf{J}_2}(m)m}\\=\quad &r_{\mathbf{J}_1(m_1)m_1}\cdot...(\operatorname{det}S^{\mathbf{J}_1}r_{\mathbf{J}_1(m_b)m_b}+\operatorname{det}S^{\mathbf{J}_2}r_{\mathbf{J}_2(m_b)m_b})...\cdot r_{\mathbf{J}_1(m_m)m_m}
\end{split}
\end{equation}
If we further assume that $\mathbf{J}_1$ and $\mathbf{J}_2$ are such that one well-behaves and the other ill-behaves we have:
\begin{equation}
\operatorname{det}S^{\mathbf{J}_1}r_{\mathbf{J}_1(m_b)m_b}+\operatorname{det}S^{\mathbf{J}_2}r_{\mathbf{J}_2(m_b)m_b}=\mathpzc{a} \cdot r_{\mathbf{J}_1(m_b)m_b}-\mathpzc{b} \cdot r_{\mathbf{J}_2(m_b)m_b},
\end{equation}
with $\mathpzc{a}$ and $\mathpzc{b}$ constants of the same sign.\\

By the mere fact that $d$ is an integer-valued distance, any other Child Selection satisfies
\begin{equation}
d(\mathbf{J}_k,\mathbf{J}_1), d(\mathbf{J}_k,\mathbf{J}_2) \ge 1,\quad \text{for any } k\neq1,2.
\end{equation}

In particular, we have the following Lemma:
\begin{lem}\label{propSNB}
Let $\mathbf{J}_1$ and $\mathbf{J}_2$ be Child Selections at distance $d=1$, that is, $\mathbf{J}_1(m_b)\neq \mathbf{J}_2(m_b)$ and $\mathbf{J}_1(m)= \mathbf{J}_2(m)$ for any $m\neq m_b$. For any other Child Selection $\mathbf{J}_k$, there is a metabolite $m_k$ such that $\mathbf{J}_k(m_k)\neq \mathbf{J}_1(m_k)$ and $\mathbf{J}_k(m_k)\neq \mathbf{J}_2(m_k)$.\\ 
Moreover if $d(\mathbf{J}_k,\mathbf{J}_1) = d(\mathbf{J}_k,\mathbf{J}_2) = 1$, then $m_k=m_b$.\\
\end{lem}
\proof
Let us consider any $m_k$ such that $\mathbf{J}_1(m_k)\neq \mathbf{J}_k(m_k)$. If $\mathbf{J}_2(m_k)\neq \mathbf{J}_k(m_k)$, we are done. Assume then that $\mathbf{J}_2(m_k)= \mathbf{J}_k(m_k)$. By construction, $m_k=m_b$. Consider now any $\tilde{m}_k$ such that $\mathbf{J}_2(\tilde{m}_k)\neq \mathbf{J}_k(\tilde{m}_k)$ and remember that $\mathbf{J}_1(m) = \mathbf{J}_2(m)$ for any $m\neq m_b$. We conclude that $\mathbf{J}_1(\tilde{m}_k)\neq \mathbf{J}_k(\tilde{m}_k)$. Otherwise, we would have found two metabolites $m_k$ and $\tilde{m}_k$ such that $\mathbf{J}_1(m_k)\neq \mathbf{J}_2(m_k)$ and $\mathbf{J}_1(\tilde{m}_k)\neq \mathbf{J}_2(\tilde{m}_k)$, contradicting $d(\mathbf{J}_1,\mathbf{J}_2)=1$.\\
In the above argument, note that if $\mathbf{J}_2(m_k)= \mathbf{J}_k(m_k)$, then $d(\mathbf{J}_1,\mathbf{J}_k)\ge2$. Hence, if $d(\mathbf{J}_1,\mathbf{J}_k) = d(\mathbf{J}_2,\mathbf{J}_k) = 1$ we conclude that $\mathbf{J}_1(m_b)\neq \mathbf{J}_2(m_b)\neq \mathbf{J}_k(m_b)$ and hence $m_k=m_b$.
\endproof

We are now ready to prove Theorem \ref{CoS}.

\proof[Proof of Theorem \ref{CoS}]
Let $\mathbf{J}_1$ and $\mathbf{J}_2$ be Child Selections, as above. In particular, distance $d(\mathbf{J}_1,\mathbf{J}_2)=1$.
By Lemma \ref{propSNB}, for any other Child Selection $\mathbf{J}_k \neq \mathbf{J}_1, \mathbf{J}_2$ we can find $m_k$ such that 
\begin{equation}
\mathbf{J}_1(m_k), \mathbf{J}_2(m_k)\neq \mathbf{J}_k(m_k).
\end{equation}
We can consider, then, an $\epsilon$-small choice of reaction rate parameters such that
\begin{equation}
r_{\mathbf{J}_k(m_k)m_k}  < \epsilon.
\end{equation}
Then, for this $\epsilon$-choice of reaction rates,
\begin{equation}
\operatorname{det}G=(\mathpzc{a}  r_{\mathbf{J}_1(m_b)m_b}-\mathpzc{b} r_{\mathbf{J}_2(m_b)m_b}) r_{\mathbf{J}_1(m_{n+1})m_{n+1}} ... r_{\mathbf{J}_1(m_m)m_m}+ \epsilon.
\end{equation}
The parameter $\xi=\mathpzc{a}  r_{\mathbf{J}_1(m_b)m_b}-\mathpzc{b} r_{\mathbf{J}_2(m_b)m_b}$ becomes then a bifurcation parameter for the sign of the Jacobian determinant.
\endproof

\bibliography{bibliography/references.bib}

\begin{thebibliography}{10}

\bibitem{HFJ72}
F.~Horn and R.~Jackson, ``General mass action kinetics,'' {\em Archive for
  Rational Mechanics and Analysis}, vol.~47, no.~2, pp.~81--116, 1972.

\bibitem{Fei87}
M.~Feinberg, ``Chemical reaction network structure and the stability of complex
  isothermal reactors—{I}. {T}he deficiency zero and deficiency one
  theorems,'' {\em Chemical Engineering Science}, vol.~42, no.~10,
  pp.~2229--2268, 1987.

\bibitem{Fei95}
M.~Feinberg, ``The existence and uniqueness of steady states for a class of
  chemical reaction networks,'' {\em Archive for Rational Mechanics and
  Analysis}, vol.~132, no.~4, pp.~311--370, 1995.

\bibitem{GaNi65}
D.~Gale and H.~Nikaido, ``The {J}acobian matrix and global univalence of
  mappings,'' {\em Mathematische Annalen}, vol.~159, no.~2, pp.~81--93, 1965.

\bibitem{CraFei06}
G.~Craciun and M.~Feinberg, ``Multiple equilibria in complex chemical reaction
  networks: {I}{I}. {T}he species-reaction graph,'' {\em SIAM Journal on
  Applied Mathematics}, vol.~66, no.~4, pp.~1321--1338, 2006.

\bibitem{Ba-07}
M.~Banaji, P.~Donnell, and S.~Baigent, ``P matrix properties, injectivity, and
  stability in chemical reaction systems,'' {\em SIAM Journal on Applied
  Mathematics}, vol.~67, no.~6, pp.~1523--1547, 2007.

\bibitem{BaCra10}
M.~Banaji and G.~Craciun, ``Graph-theoretic criteria for injectivity and unique
  equilibria in general chemical reaction systems,'' {\em Advances in Applied
  Mathematics}, vol.~44, no.~2, pp.~168--184, 2010.

\bibitem{ShiFei12}
G.~Shinar and M.~Feinberg, ``Concordant chemical reaction networks,'' {\em
  Mathematical Biosciences}, vol.~240, no.~2, pp.~92--113, 2012.

\bibitem{ShiFei13}
G.~Shinar and M.~Feinberg, ``Concordant chemical reaction networks and the
  species-reaction graph,'' {\em Mathematical Biosciences}, vol.~241, no.~1,
  pp.~1--23, 2013.

\bibitem{Thom81}
R.~Thomas, ``On the relation between the logical structure of systems and their
  ability to generate multiple steady states or sustained oscillations,'' {\em
  Numerical methods in the study of critical phenomena}, pp.~180--193, 1981.

\bibitem{Soule2003}
C.~Soul{\'e}, ``Graphic requirements for multistationarity,'' {\em ComPlexUs},
  vol.~1, no.~3, pp.~123--133, 2003.

\bibitem{KST07}
M.~Kaufman, C.~Soul{\'e}, and R.~Thomas, ``A new necessary condition on
  interaction graphs for multistationarity,'' {\em Journal of theoretical
  biology}, vol.~248, no.~4, pp.~675--685, 2007.

\bibitem{WF13}
C.~Wiuf and E.~Feliu, ``Power-law kinetics and determinant criteria for the
  preclusion of multistationarity in networks of interacting species,'' {\em
  SIAM Journal on Applied Dynamical Systems}, vol.~12, no.~4, pp.~1685--1721,
  2013.

\bibitem{Timo20}
E.~Feliu, N.~Kaihnsa, T.~de~Wolff, and O.~Y{\"u}r{\"u}k, ``The kinetic space of
  multistationarity in dual phosphorylation,'' {\em arXiv preprint
  arXiv:2001.08285}, 2020.

\bibitem{MR07}
M.~Mincheva and M.~R. Roussel, ``Graph-theoretic methods for the analysis of
  chemical and biochemical networks. {I}. {M}ultistability and oscillations in
  ordinary differential equation models,'' {\em Journal of mathematical
  biology}, vol.~55, no.~1, pp.~61--86, 2007.

\bibitem{Iv79}
A.~Ivanova and B.~Tarnopolskii, ``One approach to the determination of a number
  of qualitative features in the behavior of kinetic systems, and realization
  of this approach in a computer (critical conditions, autooscillations),''
  {\em Kinetics and Catalysis}, vol.~20, no.~6, pp.~1271--1277, 1979.

\bibitem{Iv792}
A.~Ivanova, ``Conditions for uniqueness of the stationary states of kinetic
  systems, connected with the structure of their reaction-mechanism. 1.,'' {\em
  Kinetics and Catalysis}, vol.~20, no.~4, pp.~833--837, 1979.

\bibitem{VIv87}
A.~Volpert and A.~Ivanova, ``Mathematical models in chemical kinetics,'' {\em
  Mathematical modeling (Russian)}, vol.~57, p.~102, 1987.

\bibitem{BF18}
B.~Brehm and B.~Fiedler, ``Sensitivity of chemical reaction networks: a
  structural approach. 3. {R}egular multimolecular systems,'' {\em Mathematical
  Methods in the Applied Sciences}, vol.~41, no.~4, pp.~1344--1376, 2018.

\bibitem{Fei19}
M.~Feinberg, {\em Foundations of Chemical Reaction Network Theory}.
\newblock Springer, 2019.

\bibitem{F19}
B.~Fiedler, ``Global {H}opf bifurcation in networks with fast feedback
  cycles,'' {\em Discrete and Continuous Dynamical Systems - S}, vol.~0,
  no.~1937-1632\_2019\_0\_144, 2020.

\bibitem{Vas20}
N.~Vassena, {\em Sensitivity of Metabolic Networks}.
\newblock PhD thesis, Freie Universit{\"a}t in Berlin, 2020.

\bibitem{Al03}
B.~Alberts, D.~Bray, J.~Lewis, M.~Raff, K.~Roberts, and J.~D. Watson,
  ``Molecular biology of the cell,'' 2003.

\bibitem{Lod08}
H.~Lodish, A.~Berk, C.~A. Kaiser, M.~Krieger, M.~P. Scott, A.~Bretscher,
  H.~Ploegh, P.~Matsudaira, {\em et~al.}, {\em Molecular cell biology}.
\newblock Macmillan, 2008.

\bibitem{Ish07}
N.~Ishii, K.~Nakahigashi, T.~Baba, M.~Robert, T.~Soga, A.~Kanai, T.~Hirasawa,
  M.~Naba, K.~Hirai, A.~Hoque, {\em et~al.}, ``Multiple high-throughput
  analyses monitor the response of e. coli to perturbations,'' {\em Science},
  vol.~316, no.~5824, pp.~593--597, 2007.

\bibitem{Nak+09}
K.~Nakahigashi, Y.~Toya, N.~Ishii, T.~Soga, M.~Hasegawa, H.~Watanabe, Y.~Takai,
  M.~Honma, H.~Mori, and M.~Tomita, ``Systematic phenome analysis of
  escherichia coli multiple-knockout mutants reveals hidden reactions in
  central carbon metabolism,'' {\em Molecular Systems Biology}, vol.~5, no.~1,
  p.~306, 2009.

\bibitem{VL06}
F.~G. Vital-Lopez, C.~D. Maranas, and A.~Armaou, ``Bifurcation analysis of the
  metabolism of {E}. coli at optimal enzyme levels,'' in {\em 2006 American
  Control Conference}, pp.~6--pp, IEEE, 2006.

\bibitem{Chass2002}
C.~Chassagnole, N.~Noisommit-Rizzi, J.~W. Schmid, K.~Mauch, and M.~Reuss,
  ``Dynamic modeling of the central carbon metabolism of {E}scherichia coli,''
  {\em Biotechnology and bioengineering}, vol.~79, no.~1, pp.~53--73, 2002.

\bibitem{OTM18}
T.~Okada, J.-C. Tsai, and A.~Mochizuki, ``Structural bifurcation analysis in
  chemical reaction networks,'' {\em Physical Review E}, vol.~98, no.~1,
  p.~012417, 2018.

\end{thebibliography}
 \bibliographystyle{ieeetr}

\end{document}